# On *d*-Walk Regular Graphs

Ernesto Estrada[1,2] and José A. de la Peña[2]

[1]Department of Mathematics and Statistics, University of Strathclyde, Glasgow G1 1XH, U.K.,
[2]CIMAT, Guanajuato, Mexico

ABSTRACT: Let $G$ be a graph with set of vertices $1,\ldots n$ and adjacency matrix $A = A(G) = (a_{ij})$ of size $n \times n$. Let $d(i,j) = d$, we say that $f_d : \mathbb{N} \to \mathbb{N}$ is a $d$-function on $G$ if for every pair of vertices $i, j$ and $k \geq d$, we have $a_{i,j}^{(k)} = f_d(k)$. If this function $f_d$ exists on $G$ we say that $G$ is $d$-walk regular. We prove that $G$ is $d$-walk regular if and only if for every pair of vertices $i, j$ at distance $\leq d$ and for $d \leq k \leq n + d - 1$, we have that $a_{i,j}^{(k)}$ is independent of the pair $i, j$. Equivalently, the single condition $e^A \circ A_d = cA_d$ holds for some constant $c$, where $A_d$ is the adjacency matrix of the $d$-distance graph and $\circ$ denotes the Schur product.

## 1. Introduction

Let $G$ be a graph with set of vertices $1,\ldots,n$ and adjacency matrix $A = A(G) = (a_{ij})$ of size $n \times n$. We denote by $A = \left(a_{ij}^{(k)}\right)$ the $k$th power of $A$, where $a_{ij}^{(k)}$ is the number of walks of length $k$ in $G$. Let the distance between a pair of nodes be denoted by $d(i, j) = d$. We say that $f_d : \mathbb{N} \to \mathbb{N}$ is a $d$-function on $G$ if for every pair of vertices $i, j$ and for any $k \geq d$, we have $a_{i,j}^{(k)} = f_d(k)$. If this function $f_d$ exists on $G$ we say that $G$ is $d$-walk regular.

Observe that our definition slightly differs from that given by Dalfó, Fiol and Garriga [7], where they require the condition $a_{i,j}^{(k)} = f_d(k)$ for pairs of vertices at distance $d(i, j) \leq d$. With this definition $d$-walk regular graphs have been considered in very few occasions, but other related concepts have appeared long before in the specialized literature. For instance, a graph is said to be *distance regular* if for every $k$ the virtual power $a_{ij}^{(k)}$ is constant for pairs $i, j$ at constant distance $d(i, j)$. This class of graphs includes Hamming graphs, the Johnson graphs, the

Grassmann graphs and graphs associated to dual polar spaces (see, for instance, Brouwer, Cohen and Neumaier [2]). We also mention here the Hoffman polynomial $H$ charaterizing the regularity of $G$ by the condition $H(A) = J$, where $J$ is the all-ones matrix, see [11]. Clearly, distance regular graphs are those which are $d$-walk regular for every $d$.

Given $d$ a positive integer, we consider the $d$-distance graph $G_d$ of $G$ to be the graph with vertices $1, \ldots, n$ and edges $i, j$ in case $d(i, j) = d$. We write $A_d = A(G_d) = (a_{ij}(d))$ for the adjacency matrix and $m(d)$ for the number of edges of the graph $G_d$. By convention $m(0) = n$. Recall that, for a distance-regular graph $G$, there are polynomials $p_k$ satisfying $p_k(A) = A_k$ $(0 \leq k \leq d)$, where $A_k$ stands for the adjacency matrix of the $k$-distance graph $G_k$.

0-walk regular graphs were considered by Godsil [9] under the general name of *walk-regular graphs*. The *predistance polynomials* $p_k$ were used in [7] to give the following characterization: let $G$ be a connected graph with adjacency matrix $A$ having $d+1$ distinct eigenvalues, then $G$ is 0-walk-regular if and only if the matrices $p_k(A)$, $1 \leq k \leq d$, have null diagonals. As further examples, we recall that a family of graphs which are 1-walk-regular (but not $k$-walk-regular for $k > 1$) are the Cartesian products of cycles $C_m \times C_m$ with $m \geq 5$. In [4], the current authors introduced the *unfolded complete graphs* $G$, for which there exist some $N$, called the *size* of $G$, such that for every edge $p, q$ the virtual powers of the adjacency matrix satisfy

$$a_{pq}^{(2j+1)} = \frac{(N-1)^{2j+1} + 1}{N}, \tag{1}$$

for every $j \geq 1$, which is the number of odd walks in the complete graph $K_N$. The unfolded complete graphs are 1-walk regular.

## 2. Main Results

In this section we state the main results of this work which will be proved in a further section. As a first result, we prove that $d$-walk regularity is equivalent to a finite set of conditions. Namely, we prove the following result.

**Theorem 1**: Let $G$ be a graph with adjacency matrix $A$. Then the following are equivalent:

(a) $G$ is an $d$-walk regular graph;

(b) for every pair of vertices $i, j$ at distance $d(i, j) \leq d$ and for any $d \leq k \leq n+d-1$, we have that $a_{ij}^{(k)}$ is independent of the pair $i, j$.

Moreover, in that case there is a $d$-function $f_d$ on $G$ such that for $k$ any positive number it satisfies the following equality

$$f_d(k) = \frac{1}{2m(d)} tr(A_d A^k). \tag{2}$$

Denote by $A \circ B$ the *Schur product* (sometimes also called *Hadamard or entrywise product*) of the $n \times n$ matrices $A, B$. That is $(A \circ B)_{ij} = A_{ij} B_{ij}$. On the other hand, the *Kronecker product* $A \otimes B$ is defined by the $n^2 \times n^2$ matrix $(A \otimes B)_{(ij, pq)} = A_{ip} B_{jq}$ (see [3] for a general reference). Then we have the following result.

**Theorem 2**: Let $G$ be a graph with adjacency matrix $A$. The following are equivalent:

a) $G$ is a $d$-walk regular graph;

b) for any number $1 \leq k \leq n-1$, there exists a positive integer $c(k)$ such that

$$A^k \circ A_d = c(k) A_d. \tag{3}$$

The numbers $c(k)$ are determined by

$$c(k) = \frac{1}{2m(d)} sum(A^k \circ A_d), \tag{4}$$

where $sum(a_{ij}) = \sum_{i,j=1}^{n} a_{ij}$.

c) for each number $1 \leq k \leq n-1$, there exists a positive integer $c(k)$ such that

$$A^k \otimes A_d = c(k) A_d^n, \tag{5}$$

where the numbers $c(k)$ are determined by

$$c(k) = \frac{1}{2nm(d)} sum(A^k \otimes A_d) \tag{6}$$

Our last characterization shows that the set of conditions defining walk regular graphs are equivalent to a single condition on the exponential matrix. Namely, we prove the following result.

**Theorem 3**: Let $G$ be a graph with adjacency matrix $A$. Then the following are equivalent:

(a) $G$ is a $d$-walk regular graph;

(b) there exists a positive number $c$ and a non-negative matrix $C$ such that

$$e^A A_d = cA_d + C \tag{7}$$

where the entries of $C$ are $c_{ij}(k) = 0$ whenever $d(i, j) = d$. Equivalently,

$$e^A \circ A_d = cA_d \tag{8}$$

## 3. Elementary properties and examples

We start by observing that $G$ is $d$-*walk regular* if and only if for every $k \geq d$ there is a constant $f_d(k)$ such that

$$A^k = f_d(k)A_d + B_d(k) \tag{9}$$

for a non-negative matrix $B_d(k) = \left(b_{ij}(d)^{(k)}\right)$ satisfying $b_{ij}(d)^{(k)} = 0$ whenever $d(i, j) = d$. That is,

$$b_{ij}(d)^{(k)} = \begin{cases} 0 & \text{if } a_{ij}(d) \neq 0 \\ a_{ij}^{(k)} & \text{if } a_{ij}(d) = 0. \end{cases} \tag{10}$$

Observe that for $k < d$, necessarily $f_d(k) = 0$.

Let $p(T) = T^n + p_{n-1}T^{n-1} + \cdots + p_0$ be the characteristic polynomial of the graph $G$. The Cayley-Hamilton theorem yields $p(A) = 0$. For $d \leq k_0$ we use the Euclidean algorithm to write

$$T^k = s(T)p(T) + r_k(T), \tag{11}$$

with

$$r_k(T) = r_{k,n-1}T^{n-1} + \cdots + r_{k,0}, \tag{12}$$

which is a polynomial of degree $\leq n-1$. Then

$$A^k = s(A)p(A) + r_k(A) = r_k(A) = r_{k,n-1}A^{n-1} + \cdots + r_{k,0}. \tag{13}$$

Assume that $G$ is $d$-*walk regular* and for every $0 \leq k \leq n-1$ there is a constant $c(k)$ such that

$$A^k = c(k)A_d + B_d(k), \tag{14}$$

for a non-negative matrix $B_d(k)$. Then

$$A^k = r_k(A) = (r_{k,n-1}c(n-1) + \cdots + r_{k,0}c(0))A_d - (r_{k,n-1}B_d(n-1) + \cdots + r_{k,0}B_d(0)). \tag{15}$$

That is,

$$c(k) = r_{k,n-1}c(n-1) + \cdots + r_{k,0}c(0). \tag{16}$$

Some further examples are the following.

(1) Let $S_0(N) = O_1(N), O_3(N), \ldots, O_{2j+1}(N), \ldots$ be the sequence of walks of odd length in a complete graph $K_N$. Thus

$$O_{2j+1}(N) = \frac{(N-1)^{2j+1} + 1}{N}. \tag{17}$$

Let $G$ be a graph with adjacency matrix $A$. As referred in the Introduction, we say that $G$ is *unfolded complete* if there exist some $N$, called the *size* of $G$, such that for every edge $p - q$ the virtual power $a_{pq}^{(2j+1)} = O_{2j+1}(N)$ for every $j \geq 1$. Recently, the authors proved the following result [4]. Let $G$ be a connected unfolded complete graph of size $N$. Then $G = K_2 \otimes K_N$.

Some further examples of $k$-walk regularity in graphs are provided in Figure 1. In Fig. 1A we illustrate a graph which is $k$-walk regular for $k = 0,2,4$ and it is not $k$-walk for $k = 1,3$. The graph in Fig. 1B is $k$-walk regular for $k = 0,2$ but it is not $k$-walk for $k = 1,3,4,5$.

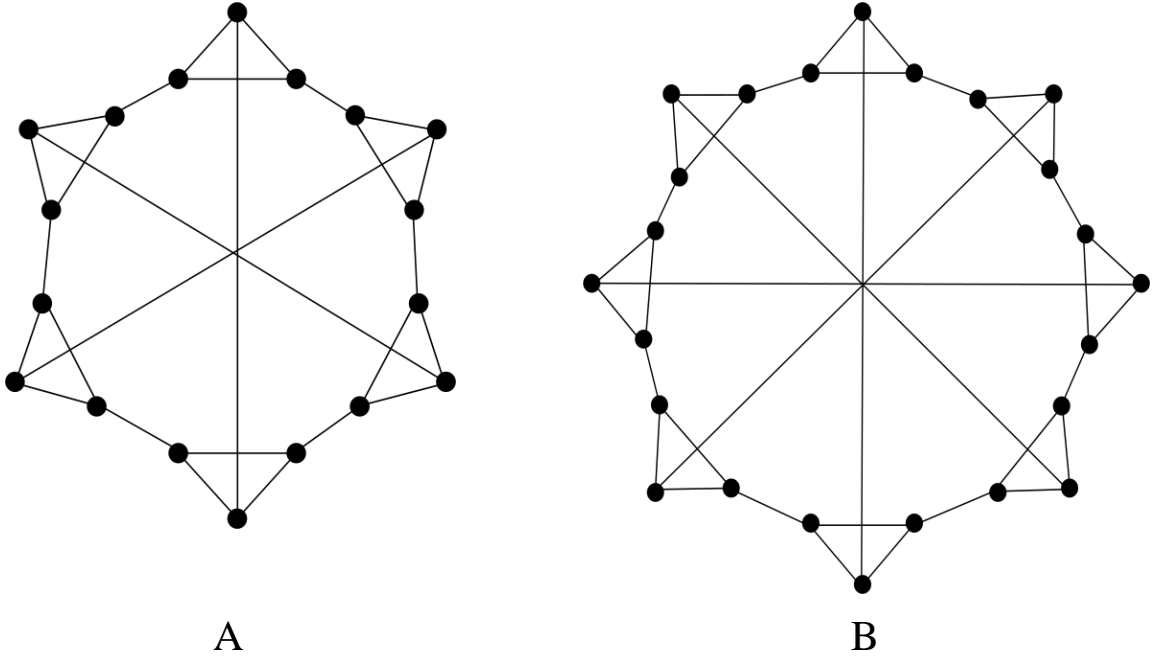

A  B

### 4. Proof of the theorems

Before starting with the proofs we recall that there are excellent books that may provide background information on the spectral graph theory techniques used in this section and the reader is referred to them for general references [1, 3, 9].

**Proof of Theorem 1**: Assume that $a_{ij}^{(k)} = f(k)$ for every pair of vertices $i, j$ at distance $d$, and let $d \leq k \leq n+1$. For those numbers $k$ let us write

$$A^k = a_{ij}^{(k)} = f(k)A_d + B(k), \qquad (18)$$

where $B(k) = \left(b_{ij}^{(k)}\right)$ is a non-negative matrix with entries defined by

$$b_{ij}^{(k)} = \begin{cases} 0 & \text{if } a_{ij}(d) \neq 0, \\ a_{ij}^{(k)} & \text{if } a_{ij}(d) = 0. \end{cases}$$

Observe that

$$(A_d B(k))_{ii} = \sum_{d(i,p)=d} a_{ip}^{(d)} b_{pi}^{(k)} = 0, \qquad (19)$$

Therefore $tr(A_d B(k)) = 0$. Moreover

$$A_d A^k = f(k) A_d^2 + A_d B(k), \qquad (20)$$

with $tr(A_d A^k) = f(k) tr(A_d^2) = 2m(d) f(k)$. The expression for $f$ follows.

We shall say that $k$ is *good* if $A^k = f(k) A_d + B(k)$ for $B(k) = (b_{ij}^{(k)})$ a non-negative matrix with entries $b_{ij}^{(k)} = 0$ if $a_{ij}(d) \neq 0$. By hypothesis, all numbers $n + d - 1 \leq k \leq N - 1$ are good. By induction we show that $N$ is good.

Let $p(T) = T^n + p_{n-1} T^{n-1} + \cdots + p_0$ be the characteristic polynomial of $A$, then the Cayley-Hamilton theorem yields, for $N \geq n$

$$0 = p(A) A^{N-n} = A^N + p_{n-1} A^{N-1} + \cdots + p_0 A^{N-n}. \qquad (21)$$

In particular,

$$f(N) + p_{n-1} f(N-1) + \cdots + p_1 f(N-n+1) + p_0 f(N-1) = 0.$$

Since the numbers $k \leq N - 1$ are good, we have that

$$\begin{aligned} A^N &= -\left(p_{n-1} f(N-1) + \cdots + p_0 f(N-n)\right) A_d - \left(p_{n-1} B(N-1) + \cdots + p_0 B(N-n)\right) \\ &= f(N) A_d + B(N), \end{aligned} \qquad (22)$$

with the matrix $B(N)$ satisfying that for each pair $i, j$ at distance $d$ we have $(B(N))_{ij} = 0$. Moreover, observe that $B(N)$ is non-negative. Therefore $N$ is good. It follows that $G$ is $d$-walk regular.

Finally, to show that $f_d$ is a $d$-function on $G$ we need to observe that it takes positive integer values. Indeed, for any $k$, if $a_{ij}(d) \neq 0$ then $a_{ij}^{(k)} = f(k) a_{ij}(d) = f(k)$ is a positive integer.

□

***Proof of Theorem 2:*** We show that (a) implies (b) in Theorem 2: suppose that $A^k = f(k)A_d + B(k)$, where $B(k) = \left(b_{ij}^{(k)}\right)$ is a non-negative matrix with entries $b_{ij}^{(k)} = 0$ if $a_{ij}(d) \neq 0$. Then, $A_d \circ A_d = A_d$ and $A_d \circ B(k) = 0$. Therefore $A^k \circ A_d = f(k)A_d$. Assume that (b) holds. Take a number $1 \leq k \leq n-1$, and write

$$A^k = c(k)A_d + B_d(k),$$

for $B(k) = \left(b_{ij}^{(k)}\right)$ an $n \times n$ matrix. Since $A^k \circ A_d = c(k)A_d$ then $B(k) \circ A_d = 0$. That is, for every edge $(i, j) \in E$ in $G_d$ we have $b_{ij}^{(k)} = 0$. Hence (a) holds.

We now show that (a) implies (d): Suppose that $A^k = f(k)A_d + B(k)$ for $B(k) = \left(b_{ij}^{(k)}\right)$ a non-negative matrix with entries $b_{ij}^{(k)} = 0$ if $a_{ij}(d) \neq 0$. Then, $A_d \otimes A_d = A_d^{\ n}$ and $A_d \otimes B(k) = f(k)A_d^{\ n}$ are matrices of size $n^2 \times n^2$. Assume (d) holds. Take a number $1 \leq k \leq n-1$, and write

$$A^k = c(k)A_d + B(k),$$

for $B(k) = \left(b_{ij}^{(k)}\right)$ an $n \times n$ matrix. Since $A^k \otimes A_d = c(k)A_d^{\ n}$ then $B(k) \otimes A_d = 0$. That is, for every edge $(i, j) \in E$ in $G_d$ we have $b_{ij}^{(k)} = 0$. Hence (a) holds. □

***Proof of Theorem 3:*** We show that (a) implies (b): Suppose that $A^k = f(k)A_d + B(k)$, where $B(k) = \left(b_{ij}^{(k)}\right)$ is a non-negative matrix with entries $b_{ij}^{(k)} = 0$ if $a_{ij}(d) \neq 0$. Then $A^k A^d = f(k)A_d + C(k)$, where $C(k)$ is a matrix with entries $c_{ij}(k) = 0$ whenever $d(i, j) = d$, for all $k \geq 0$. Hence,

$$e^A A_d = \sum_{k=0}^{\infty} \frac{A^k A_d}{k!} + \left(\sum_{k=0}^{\infty} \frac{C(k)}{k!}\right) = \left(\sum_{k=0}^{\infty} \frac{f(k)}{k!}\right) A_d + C, \qquad (23)$$

where the entries of $C$ are $c_{ij} = 0$ whenever $d(i, j) = d$.

Assume (b) holds, that is, $e^A A^d = cA_d + C$, with the entries of the non-negative matrix $C$ satisfying the desired condition. Take a number $0 \leq k$, and write

$$A^k = c(k)A_d + B(k), \qquad (24)$$

for some $c(k) \geq 0$ maximal such that $B(k) = \left(b_{ij}^{(k)}\right)$ is a non-negative matrix. Then

$$cA_d + C = e^A A_d = \left(\sum_{k=0}^{\infty} \frac{c(k)}{k!}\right) A_d + \left(\sum_{k=0}^{\infty} \frac{B(k)}{k!}\right). \qquad (25)$$

Thus,

$$\left(\sum_{k=0}^{\infty} \frac{c(k)}{k!}\right) \leq c, \qquad (26)$$

and if equality does not hold, then there is some $k$ and $\varepsilon > 0$ such that $b_{ij}^{(k)} \geq \varepsilon a_{ij}(d)$, contradicting the choice of $c(k)$. Consequently, the equality holds and

$$C = \sum_{k=0}^{\infty} \frac{B(k)}{k!}. \qquad (27)$$

Since $B(k)$ is non-negative, then for every edge $(i, j) \in E$ in $G_d$ we have $b_{ij}^{(k)} = 0$. Hence (a) holds. □

### 5. Additional remarks

**Proposition 4**: (1) Assume that condition (b) in Theorem 2 holds. Then the following two equivalent assertions are true:

(a) for each number $d \leq k \leq n + d - 1$, there exists a positive integer $c(k)$ such that any non-zero eigenvalue $\vartheta$ of the Schur product $A^k \circ A_d$ is of the form $\vartheta = \mu c(k)$ for some eigenvalue $\mu$ of $A_d$.

(b) for each number $d \leq k \leq n+d-1$, there exists a positive integer $c(k)$ such that the characteristic polynomial of the Schur matrix $A^k \circ A_d$ is of the form

$$p_{A^k \circ A_d}(x) = \sum_{i=0}^{n} b_i x^i \text{ with } b_i = a_i c(k)^i \tag{28}$$

for every $1 \leq i \leq n$, where $p_{A_d}(x) = \sum_{i=0}^{n} a_i x^i$.

(2) Assume that assertion (d) holds, the following two equivalent assertions are true:

(c) for each number $d \leq k \leq n+d-1$, there exists a positive integer $c(k)$ such that any non-zero eigenvalue $\vartheta$ of the tensor product $A^k \otimes A_d$ is of the form $\vartheta = \mu c(k)$ for some eigenvalue $\mu$ of $A_d$.

(d) for each number $d \leq k \leq n+d-1$, there exists a positive integer $c(k)$ such that the characteristic polynomial of the tensor product $A^k \otimes A_d$ is of the form

$$p_{A^k \otimes A_d}(x) = \sum_{i=0}^{n} b_i x^i \text{ with } b_i = a_i c(k)^i \tag{29}$$

for every $1 \leq i \leq n$, where $p_{A_d}(x) = \sum_{i=0}^{n} a_i x^i$.

*Proof:* We show first that (b) in Theorem 2 implies (a) and (b) in Proposition 4. Choose an eigenvector $u$ of $A_d$ with eigenvalue $\mu \neq 0$. That is, $A_d u = \mu u$ and

$$(A^k \circ A_d) u = \mu f(k) u, \tag{30}$$

Hence, if $p_{A_d}(x) = \prod_{\mu}(x - \mu^{e(\mu)}) = \sum_{i=0}^{n} a_i x^i$ then

$$p_{A^k \circ A_d}(x) = \prod_{\mu}(x - \mu^{e(\mu)} f(k)) = \sum_{i=0}^{n} b_i x^i. \tag{31}$$

Therefore $b_i = a_i f(k)^i$ for every $1 \leq i \leq n$.

Since the coefficient $a_i$ (resp. $b_i$) is the $i$-th symmetric function on the eigenvalues of $A_d$ (resp. of $A^k \circ A_d$), then (a) implies (b) in Proposition 4. It remains to show that (b) implies (a) in Proposition 4. Indeed, assume that (b) holds for $2 \leq k \leq n-1$, then

$$p_{A^k \circ A_d}\left(c(k)^{-1} A_d\right) = \sum_{i=0}^{n} b_i c(k)^{-1} A_d = \sum_{i=0}^{n} a_i A_d^i = p_{A_d}(A_d) = 0. \tag{32}$$

A further application of the Cayley-Hamilton theorem assures that the minimal polynomial $q_{A^k \circ A_d}(x)$ of $A^k \circ A_d$ is a divisor of $p_{A_d}\left(c(k)^{-1} x\right)$. This yields that any non-zero eigenvalue $\vartheta$ of $A^k \circ A_d$ is of the form $\vartheta = \mu c(k)$ for some eigenvalue $\mu$ of $A_d$. This completes the proof of (1). Similarly, we may show (2). □

**6. 0-walk regular graphs and regularity.**

Let us assume that $G$ is $d$-walk regular, then for $A = A(G)$ and for any $k \geq d$ there is a constant $c(k)$ such that

$$A^k = c(k) A_d + B_d(k), \tag{33}$$

for a matrix $B(k) = \left(b_{ij}(d)^{(k)}\right)$ satisfying $b_{ij}(d)^{(k)} = 0$ whenever $d(i,j) = d$. Observe that for $k < d$ the above equation holds for $c(k) = 0$. Denote by $N_i(d)$ the number of vertices at distance $d$ from $i$, that is, the number of neighbours of $i$ in the graph $G_d$.

**Proposition 5**: Let $G$ be a $d$-walk regular graph. Then the following holds:

(a) If $d = 0$ then $G$ is a regular graph.
(b) If $d > 0$ then $G$ is a 0-walk regular graph if and only if the graph $G_d$ is regular.

*Proof*: For (a), assume that $d = 0$. Then there is a constant $c$ such that

$$\sum_{i \sim j} a_{ij} = a_{ii}^{(2)} = c, \tag{34}$$

for every vertex $i$.

For (b), assume that for every $k \geq d$ there is a constant $f_d(k)$ such that

$$A^k = c(k)A_d + B_d(k), \qquad (35)$$

for a matrix $B(k) = (b_{ij}(d)^{(k)})$ satisfying $b_{ij}(d)^{(k)} = 0$ whenever $d(i,j) = d$. Calculating the $i$th diagonal entry of the product of the above equation and $A^k$ (respectively $A^{k+1}$) yields

$$a_{ii}^{(2k)} = c(k)\sum_{s=1}^{n} a_{is}^{(k)} a_{is}(d) = c(k)^2 N_i(d), \qquad (36)$$

$$a_{ii}^{(2k+1)} = c(k+1)\sum_{s=1}^{n} a_{is}^{(k+1)} a_{is}(d) = c(k+1)^2 N_i(d). \qquad (37)$$

Observe that in case $N_i(d)$ is independent of the vertex $i$, then we can write

$$A^k = r(k)I_n + B(k), \qquad (38)$$

where $B(k)$ is a hollow matrix. Conversely, if $a_{ii}^{(k)}$ depends only on $k$ then $N_i(d)$ is independent of the chosen vertex. □

Now, let us denote by $p_G(x)$ the characteristic polynomial of $G$. For a vertex $i$ denote $A_i$ the adjacency matrix of the graph obtained by deleting the row and column corresponding to $i$. Let $S$ be any set of walks in $G$. For the number of walks in $S$ of length $r$ write $S^{(r)}$ and $S(x) = \sum_r S^r x^r$ the corresponding generating function. For instance, if $C_{G,i}$ is the set of closed walks starting (and ending) at vertex $i$, then $G$ is 0-walk regular if and only if $C_{G,i}(x) = C_{G,j}(x)$ for every $i, j$. Godsil and McKay [10] have shown that $G$ is 0-walk regular if and only if the characteristic polynomial $p_{G_i}(x)$ of the graph $G_i = G \setminus \{i\}$ is independent of the chosen vertex $i$. Indeed, this follows from the formula:

$$p_{G_i}(x) = \frac{1}{x} p_G(x) C_{G,i}\left(\frac{1}{x}\right). \qquad (39)$$

We recall here that the *line graph* $L(G)$ of $G$ is the graph whose vertices are the edges of $G$ and there is an edge $p,q$ whenever $p$ and $q$ have a common vertex. Observe that a walk $i_1 - i_2 - \cdots i_s - i_{s+1}$ of length $s$ in $G$ yields a path $(i_1-i_2)-(i_2-i_3)\cdots(i_s-i_{s+1})$ of length $s-1$ in $L(G)$.

In what follows we denote by $G_i = G \setminus \alpha$ the graph obtained by deleting an edge $\alpha$.

**Proposition 6:** The following are equivalent for a regular graph $G$:

(a) $G$ is 1-walk regular;
(b) $L(G)$ is 0-walk regular;
(c) The characteristic polynomial $p_{L(G\setminus\alpha)}(x)$ is independent of the deleted edge $\alpha$.

*Proof:* The equivalence of (a) and (b) is clear.

Assume that $G$ is 1-walk regular. By [10], $p_{L(G\setminus\alpha)}(x)$ is the same polynomial independently of the deleted vertex $\alpha$ of $L(G)$. Let $\alpha$ be an edge in $G$ then since $L(G\setminus\alpha) = L(G)\setminus\{\alpha\}$ the statement (c) follows. Similarly, (c) implies (a). □

### Acknowledgements

This work was done during a visit of EE to CIMAT at Guanajuato. Both authors acknowledge support of CIMAT.